\documentclass[12pt]{article}
\usepackage{mathrsfs}
\usepackage{epic,eepic,epsf,epsfig}
\usepackage{amsfonts,srcltx,mathrsfs}
\usepackage{amsmath}
\usepackage{amssymb}
\usepackage{amsbsy}
\usepackage{float}
\usepackage{graphicx}
\usepackage{amsfonts}
\usepackage{color}

\usepackage{url}
\newtheorem{lem}{Lemma}[section]
\newtheorem{theorem}[lem]{Theorem}

\newtheorem{cor}[lem]{Corollary}

\newtheorem{prob}[lem]{Problem}

\newtheorem{case}{Case}
\newtheorem{prop}[lem]{Proposition}

\def\r{\rho} \def\s{\sigma}   
 \def\ld{\lambda} 
 \def\Ga{\Gamma}

 \def\lg{\langle} \def\rg{\rangle}
\def\nd{\mathrel{\bigm|\kern-.7em/}}

\def\f{\noindent}

\def\PSL{\hbox{\rm PSL}}

\def\Aut{\hbox{\rm Aut}}

\def\Aut{\hbox{\rm Aut}}

\def\demo{\f {\bf Proof.}\hskip10pt}

\newcommand{\qed}{\mbox{\raisebox{0.7ex}{\fbox{}}} \vspace{4truemm}}

\def\P{\mathcal {P}}
\def\Q{\mathcal {Q}}

\begin{document}

\title{Two infinite families of chiral polytopes of type \{4,4,4\} with solvable automorphism groups}

\author{ \\Marston D.E. Conder$^{a}$, Yan-Quan Feng$^{b}$, Dong-Dong Hou$^{b}$\\[+4pt]
$^{a}${\small Department of Mathematics, University of Auckland, PB 92019, Auckland 1142, New Zealand}\\
$^{b}${\small {Department of Mathematics, Beijing Jiaotong University, Beijing,
100044, P.R. China}}}

\date{}
\maketitle

\footnotetext{E-mails: m.conder$@$auckland.ac.nz, yqfeng$@$bjtu.edu.cn, 16118416$@$bjtu.edu.cn
}
\begin{abstract}

We construct two infinite families of locally toroidal chiral polytopes of type $\{4,4,4\}$,
with $1024m^2$ and $2048m^2$ automorphisms for every positive integer $m$, respectively.
The automorphism groups of these polytopes are solvable groups, and when $m$ is a power of $2$,  
they provide examples with automorphism groups of order $2^n$ where $n$ can be any integer greater than $9$. 
(On the other hand, no chiral polytopes of type $[4,4,4]$ exist for $n \leq 9$.) 
In particular,  our two families give a partial answer to a problem proposed by Schulte and Weiss
in [Problems on polytopes, their groups, and realizations, {\em Periodica Math.\ Hungarica\/} 53 (2006), 231-255].

\bigskip
\f {\bf Keywords:} Chiral $4$-polytope, locally toroidal, solvable group.
\end{abstract}

\section{Introduction}

Abstract polytopes are combinatorial structures with properties that generalise those of
classical polytopes.  In many ways they are more fascinating than convex polytopes
and tessellations.  Highly symmetric examples of abstract polytopes include not
only classical regular polytopes such as the Platonic solids, and more exotic structures such as
the $120$-cell and $600$-cell, but also regular maps on surfaces (such as Klein's quartic).

Roughly speaking, an abstract polytope $\P$ is a partially ordered set endowed with
a rank function, satisfying certain conditions that arise naturally from a geometric
setting. Such objects were proposed by Gr\"unbaum in the 1970s, and their definition
(initially as `incidence polytopes') and theory were developed by Danzer and Schulte.

An \emph{automorphism} of an abstract polytope $\P$ is an order-preserving permutation
of its elements, and every automorphism of is uniquely determined by its effect
on any maximal chain in $\P$ (which is known as a `flag' in $\P$). The
most symmetric examples are regular, with all flags lying in a single orbit.
The comprehensive book written by Peter McMullen and Egon Schulte~\cite{ARP} is
nowadays seen as the principal reference on this subject.

An interesting class of examples which are not quite regular are the chiral polytopes.
For these, the automorphism group has two orbits on flags, with any two flags that
differ in a single element lying in different orbits.
Chirality is a fascinating phenomenon that does not have a counterpart in the classical
theory of traditional convex polytopes.
The study of chiral abstract polytopes
was pioneered by Schulte and Weiss (see~\cite{ChiralPolytopes, ChiralPolytopesWithPSL} for example),
but it has been something of a challenge to find and construct finite examples.

For quite some time, the only known finite examples of chiral polytopes had ranks $3$ and $4$.
In rank $3$, these are are given by the irreflexible (chiral) maps on closed
compact surfaces (see Coxeter and Moser~\cite{BooksCoxeter}).
Some infinite examples of chiral polytopes of rank $5$ were constructed by Schulte
and Weiss in~\cite{FEChiralPolytopes}, and then some finite examples
of rank $5$ were constructed just over ten years ago by Conder, Hubard and Pisanski~\cite{ConstructChiralPolytopes}.

Many small examples of chiral polytopes are now known.
These include all chiral polytopes with at most $4000$ flags,
and all that are constructible from an almost simple group $\Ga$ 
of order less than $900000$.
These have been assembled in collections, as in~\cite{atles1, atles2}, for example.
In early $2009$ Conder and Devillers devised a construction for chiral polytopes whose facets are simplices,
and used this to construct examples of finite chiral polytopes of ranks $6$, $7$ and $8$ [unpublished].

 At about the same time, Pellicer devised a quite different method for constructing finite chiral polytopes
with given regular facets, and used this construction to prove the existence of finite chiral polytopes of every rank $d \geq 3$; see~\cite{ConstructChiralPolytopesforanyrank}.
A few years later, Cunningham and Pellicer proved that every finite chiral $d$-polytope with regular facets
is itself the facet of a chiral $(d +1)$-polytope; see~\cite{moreonConstructChiralPolytopesforanyrank}.

The work of Conder and Devillers was later taken up by Conder, Hubard, O'Reilly Regueiro and Pellicer~\cite{ChiralPolytopesOfANSN},
to prove that all but finitely many alternating groups $A_n$ and symmetric groups $S_n$ are the automorphism group
of a chiral 4-polytope of type $\{3,3,k\}$ for some $k$ (dependent on $n$).
This has recently been extended to every rank greater than $4$ by the same authors as~\cite{ChiralPolytopesOfANSN}.

Also Conder and Zhang in~\cite{AbelianCoverOfChiralPolytopes} introduced a new covering method
that allows the construction of some infinite families of chiral polytopes,
with each member of a family having the same rank as the original, but with the size of the members of the family
growing linearly with one (or more) of the parameters making up its `type' (Schl\"afli symbol).
They have used this method to construct several new infinite families of chiral polytopes of ranks $3, 4, 5$ and $6$.
Furthermore, Zhang constructed in her PhD thesis~\cite{Zhang'sPhdThesis} a number of chiral polytopes
of types $\{4, 4, 4\}$, $\{4, 4, 4, 4\}$ and $\{4, 4, 4, 4, 4\}$, with automorphism groups of orders $2^{10}, 2^{11}, 2^{12}$,
and $2^{15}, 2^{16}, \cdots, 2^{22}$, and $2^{18}, 2^{19}$, respectively.

\smallskip
Now let $\P$ be a regular or chiral $4$-polytope. We say that $\P$ is {\em locally toroidal} if its facets and its vertex-figures are
maps on the $2$-sphere or on the torus, and either its facets or its vertex-figures (or both) are toroidal.
Up to duality, rank 4 polytopes that are locally toroidal are of type $\{4,4,3\}$, $\{4,4,4\}$, $\{6,3,3\}$,
$\{6,3,4\}$, $\{6,3,5\}$, $\{6,3,6\}$ or $\{3,6,3\}$.

Schulte and Weiss~\cite{ChiralPolytopesWithPSL} developed a construction that starts with a $3$-dimensional
regular hyperbolic honeycomb and a faithful representation of its symmetry group as a group of complex M\"obius transformations
(generated by the inversions in four circles that cut one another at the same angles as the corresponding
reflection planes in hyperbolic space), and then derived chiral $4$-polytopes by applying
modular reduction techniques to the corresponding matrix group (see Monson and Schulte~\cite{BE2009}). 
They then used the simple group $\PSL(2, p)$ with $p$ an odd prime to construct infinite families of such polytopes.

Some years later, Breda, Jones and Schulte~\cite{ChiralPolytopesofPSLextension} developed a method of `mixing' a
chiral $d$-polytope with a regular $d$-polytope to produce a larger example of a chiral polytope of the same rank $d$.
They used this to construct such polytopes with automorphism group $\PSL(2, p) \times \Ga^{+}(\Omega)$, where $\Ga^{+}(\Omega)$
is the rotation group of a finite regular locally-toroidal $4$-polytope.
For example, $\Ga^{+}(\Omega)$ could be $A_6 \times C_2$ or $A_5$, when the corresponding
chiral polytope has type $\{4, 4, 3\}$ or $\{6, 3, 3\}$, respectively.

One can see that almost all of the examples mentioned above involve non-abelian simple groups.
On the other hand, there appear to be few known examples of chiral polytopes with solvable automorphism groups, 
apart from some of small order, and families of rank 3 polytopes arising from chiral maps on the torus (of type $\{3,6\}$, $\{4,4\}$ or $\{6,3\}$). 
This was the main motivation for the research leading to this paper.
It was also motivated in part by a problem posed by Schulte and Weiss~\cite{Problem}, namely the following:
\begin{prob}
Characterize the groups of orders $2^n$,  with $n$ a positive integer, which are automorphism groups of regular or chiral polytopes.
\end{prob}

Here we construct two infinite families of locally toroidal chiral 4-polytopes of type $\{4,4,4\}$, with solvable automorphism groups.
Each family contains one example with $1024m^2$ or $2048m^2$ automorphisms, respectively, for every integer $m \ge 1$.
In particular, if we let $m$ be an arbitrary power of $2$, say $2^k$ (with $k \ge 0$), then the automorphism group
has order $2^{10+2k}$ or $2^{11+2k}$, which can be expressed as $2^n$ for an arbitrary integer $n \ge 10$.

This extends some earlier work by the second and third authors (in~\cite{HFL} and~\cite{HFL1}), which showed that
all conceivable ranks and types can be achieved for {\em regular\/} polytopes with automorphism group of $2$-power order.
It also extends both the work by Zhang~\cite{Zhang'sPhdThesis} mentioned above, and a construction by Cunningham~\cite{TightChiralPolyhedra}
of infinite families of tight chiral $3$-polytopes of type $\{k_1, k_2\}$ with automorphism group of order $2k_1k_2$
(considered here for the special cases where $k_1$ and $k_2$ are powers of $2$).

\section{Additional background}
\label{background}

In this section we give some further background that may be helpful for the rest of the paper.

\subsection{Abstract polytopes: definition, structure and properties}

An abstract polytope of rank $n$ is a partially ordered set $\P$ endowed with a strictly monotone rank function
with range $\{-1, 0, \cdots, n\}$, which satisfies four conditions, to be given shortly.

The elements of $\P$ are called \emph{faces} of $\P$. More specifically, the elements of $\P$ of rank $j$ are called $j$-faces,
and a typical $j$-face is denoted by $F_j$.
Two faces $F$ and $G$ of $\P$ are said to be \emph{incident} with each other if $F \leq G$ or $F \geq G$ in $\P$.
A \emph{chain} of $\P$ is a totally ordered subset of $\P$, and is said to have \emph{length} $i$ if it contains exactly $i+1$ faces.
The maximal chains in $\P$ are called the \emph{flags} of $\P$.
Two flags are said to be $j$-\emph{adjacent} if they differ in just one face of rank $j$,
or simply \emph{adjacent} (to each other) if they are $j$-adjacent for some~$j$.
Also if $F$ and $G$ are faces of $\P$ with $F \leq G$, then the set $\{\,H \in \P \ |\ F \leq H \leq G\,\}$ is
called a \emph{section} of $\P$, and is denoted by $G/F$.  Such a section has rank $m-k-1$,
where $m$ and $k$ are the ranks of $G$ and $F$ respectively.
A section of rank $d$ is called a $d$-section.

\smallskip
We can now give the four conditions that are required of $\P$ to make it an abstract polytope.
These are listed as (P1) to (P4) below:
\begin{itemize}
\item [(P1)]  $\P$ contains a least face and a greatest face, denoted by $F_{-1}$ and $F_n$, respectively.
\item [(P2)]  Each flag of $\P$ has length $n+1$ (so has exactly $n+2$ faces, including $F_{-1}$ and $F_n$).
\item [(P3)]  $\P$ is \emph{strong flag-connected}, which means that any two flags $\Phi$ and $\Psi$ of $P$ can
be joined by a sequence of successively adjacent flags $\Phi=\Phi_{0}, \Phi_1, \cdots, \Phi_k=\Psi$, each of which contains $\Phi \cap \Psi$.
\item [(P4)]  The rank $1$ sections of $\P$ have a certain homogeneity property known as the \emph{diamond condition},
namely as follows: if $F$ and $G$ are incidence faces of $\P$, of ranks $i-1$ and $i+1$, respectively, where $0 \le i \le n-1$,
then there exist precisely \emph{two} $i$-faces $H$ in $\P$ such that $F< H< G$.
\end{itemize}

\noindent
An easy case of the diamond condition occurs for polytopes of rank 3 (or polyhedra): if $v$ is a vertex of same face $f$,
then there are two edges that are incident with both $v$ and $f$.

\smallskip
Next, every $2$-section $G/F$ of $\P$ is isomorphic to the face lattice of a polygon.
Now if it happens that the number of sides of every such polygon depends only on the rank of $G$,
and not on $F$ or $G$ itself, then we say that the polytope $\P$ is {\em equivelar}.
In this case, if $k_i$ is the number of edges of every $2$-section between an $(i-2)$-face and an $(i+1)$-face of $\P$,
for $1 \leq i \leq n$, then the expression $\{k_1, k_2, \cdots, k_{n-1}\}$ is called the Schl\"afli type of $\P$.
(For example, if $\P$ has rank 3, then $k_1$ and $k_2$ are the valency of each vertex and the size of each 2-face,
respectively.)

\subsection{Automorphisms of polytopes}

An \emph{automorphism} of an abstract polytope $\P$ is an order-preserving permutation of its elements.
In particular, every automorphism preserves the set of faces of any given rank.
Under permutation composition, the set of all automorphisms of $\P$ forms a group, called the automorphism
group of $\P$, and denoted by $\Aut(\P)$ or sometimes more simply as $\Ga(\P)$.
Also it is not difficult to use the diamond condition and strong flag-connectedness to prove that if an automorphism
preserves of flag of $\P$, then it fixes every flag of $\P$ and hence every element of $\P$.
It follows that $\Ga(\P)$ acts semi-regularly (or fixed-point-freely) on $\P$.

\smallskip
A polytope $\P$ is said to be \emph{regular} if its automorphism group $\Ga(\P)$ acts transitively
(and hence regularly) on the set of flags of $\P$.
In this case, the number of automorphisms of $\P$ is as large as possible, and equal to the number of flags of $\P$.
In particular, $\P$ is equivelar, and the stabiliser in $\Ga(\P)$ of every 2-section of $\P$ induces the full dihedral group
on the corresponding polygon.
Moreover, for a given flag $\Phi$ and for every $i \in \{0,1,\dots,n-1\}$, the polytope $\P$ has a unique
automorphism $\rho_i$ that takes $\Phi$ to the unique flag $(\Phi)^{i}$ that differs from $\Phi$ in precisely its $i$-face,
and then the automorphisms $\rho_0,\rho_1,\dots,\rho_{n-1}$ generate $\Ga(\P)$ and satisfy the defining
relations for the string Coxeter group $[k_1, k_2, \cdots, k_{n-1}]$, where the $k_i$ are as given in the
previous subsection for the Schl\"afli type of $\P$.
They also satisfy a certain `intersection condition', which follows from the diamond and strong flag-connectedness
conditions.
These and many more properties of regular polytopes may be found in~\cite{ARP}.

\smallskip
We now turn to chiral polytopes, for which a good reference is~\cite{ChiralPolytopes}.

\smallskip
A polytope $\P$ said to be {\em chiral\/} if its automorphism group $\Ga(\P)$ has two orbits on flags,
with every two adjacent flags lying in different orbits.
(Another way of viewing this definition is to consider $\P$ as admitting no `reflectiing' automorphism
that interchanges a flag with an adjacent flag.)
Here the number of flags of $\P$ is $2|\Ga(\P)|$,  and $\Ga(\P)$ acts regularly on each of two orbits.
Again $\P$ is equivelar, with the stabiliser in $\Ga(\P)$ of every 2-section of $\P$ inducing the full cyclic group on the corresponding polygon.

Next, for a given flag $\Phi$ and for every $j \in \{1,2,\dots,n-1\}$, the chiral polytope $\P$ admits an
automorphism $\s_j$ that takes $\Phi$ to the flag $(\Phi)^{j, j-1}$ which differs from $\Phi$ in precisely its $(j-1)$-and $j$-faces.
These automorphisms $\s_1,\s_2,\dots,\s_{n-1}$ generate $\Ga(\P)$, and if $\P$ has Schl\"afli type $\{k_1, k_2, \dots, k_{n-1}\}$,
then they satisfy the defining relations for the orientation-preserving subgroup of (index $2$ in) the
string Coxeter group $[k_1, k_2, \cdots, k_{n-1}]$.
Also they satisfy a `chiral' form of the intersection condition, which is a variant of the one mentioned earlier for regular polytopes.

Chiral polytopes occur in pairs (or {\em enantiomorphic\/} forms), such that each member of the pair is the `mirror image' of the other.
Suppose one of them is $\P$, and has Schl\"afli type $\{k_1, k_2, \cdots, k_{n-1}\}$. Then $\Ga(\P)$ is isomorphic
to the quotient of the orientation-preserving subgroup $\Lambda^{\rm o}$ of the Coxeter group $\Lambda = [k_1, k_2, \cdots, k_{n-1}]$
via some normal subgroup $K$.  By chirality, $K$ is not normal in the full Coxeter group $\Lambda$, but is conjugated by
any orientation-reversing element $c \in \Lambda$ to another normal subgroup $K^c$ which is the kernel of an epimorphism from $\Lambda^{\rm o}$
to the automorphism group $\Ga(\P^c)$ of the mirror image $\P^c$ of $\P$.

The automorphism groups of $\P$ and $\P^c$ are isomorphic to each other, but their canonical generating sets
satisfy different defining relations.   In fact, replacing the elements $\s_1$ and $\s_2$ in the canonical generating
tuple $(\s_1,\s_2,\s_3,\dots,\s_{n-1})$ by $\s_1^{-1}$ and $\s_1^{\,2}\s_2$ gives a set of generators for $\Ga(\P)$
that satisfy the same defining relations as a canonical generating tuple for $\Ga(\P^c)$,
but chirality ensures that there is no automorphism of $\Ga(\P)$ that takes $(\s_1,\s_2)$ to $(\s_1^{-1},\s_1^{\,2}\s_2)$
and fixes all the other $\s_j$.

Conversely, any finite group $G$ that is generated by $n-1$ elements $\s_1,\s_2,\dots,\s_{n-1}$ which satisfy both
the defining relations for $\Lambda^{\rm o}$ and the chiral form of the intersection condition is the `rotation subgroup' of an
abstract $n$-polytope $\P$ that is either regular or chiral.
Indeed, $\P$ is regular if and only if $G$ admits a group automorphism $\r$ of order $2$
that takes $(\s_1,\s_2,\s_3,\dots,\s_{n-1})$ to $(\s_1^{-1},\s_1^{\,2}\s_2,\s_3,\dots,\s_{n-1})$.

\smallskip
We now focus our attention on the rank 4 case.
Here the generators $\s_1, \s_2, \s_3$  for $\Ga(\P)$ satisfy the canonical relations
$\,\s_1^{k_1} = \s_2^{k_2} = \s_3^{k_3} = (\s_1\s_2)^2 = (\s_2\s_3)^2 = (\s_1\s_2\s_3)^2 = 1$,
and the chiral form of the intersection condition can be abbreviated to
$$\lg \s_1\rg \cap \lg \s_2, \s_3\rg = \{1\} = \lg \s_1, \s_2\rg \cap \lg \s_3\rg \ \hbox{ and } \ \lg \s_1, \s_2\rg \cap \lg \s_2, \s_3\rg =\lg \s_2\rg.$$

The following proposition is useful for the groups we will deal with in the proof of our main theorem.
It is called the {\em quotient criterion} for chiral $4$-polytopes.

\begin{prop}{\rm \cite[Lemma 3.2]{ChiralPolytopesofPSLextension}}\label{quotient criterion}
$\,$
Let $G$ be a group generated by elements $\s_1, \s_2, \s_3$ such that $(\s_1\s_2)^2=(\s_2\s_3)^2=(\s_1\s_2\s_3)^2=1$,
and let $\theta\!: G \to H$ be a group homomorphism taking $\,\s_j \mapsto \ld_j\,$ for $1\le j \le 3$,
such that the restriction of $\theta$ to either $\lg \s_1, \s_2\rg$ or $\lg \s_2, \s_3\rg$ is injective.
If $(\ld_{1}, \ld_{2}, \ld_{3})$ is a canonical generating triple for $H$ as the automorphism group of some chiral $4$-polytope,
then the triple $(\s_1, \s_2, \s_3)$ satisfies the chiral form of the intersection condition for $G$.
\end{prop}

\subsection{Group theory}

We use  standard notation for group theory, as in~\cite{GroupBook} for example.
In this subsection we briefly describe some of the specific aspects of group theory that we need.

\smallskip
Let $G$ be any group.   We define the {\em commutator\/} $[x, y]$ of elements $x$ and $y$ of $G$
by $[x, y]=x^{-1}y^{-1}xy$, and then define the {\em derived subgroup} (or commutator subgroup) of $G$
as the subgroup $G'$ of $G$ generated by all such commutators.
Then for any non-negative integer $n$, we define the $n\,$th derived group of $G$ by setting
$$G^{(0)}=G, \ G^{(1)}=G', \ \hbox{and } \ G^{(n)}=(G^{(n-1)})' \ \hbox{when} \ n \geq 1.$$

A group $G$ is called {\em solvable\/}  if $G^{(n)}=1$ for some $n$.
(This terminology comes from Galois theory, because a polynomial over a field $\mathbb{F}$ is solvable by radicals
if and only if its Galois group over $\mathbb{F}$ is a solvable group.)
Every abelian group and every finite $p$-group is solvable, but every non-abelian simple groups is not solvable.
In fact, the smallest non-abelian simple group $A_5$ is also the smallest non-solvable group.

We also need the following, which are elementary and so we give them without proof.

\begin{prop}\label{solvable}
If $N$ is a normal subgroup of a group $G$, such that both $N$ and $G/N$ are solvable, then so is $G$. 
\end{prop}

\begin{prop}\label{freeabelian}
Let $G$ be the free abelian group $\mathbb{Z}\oplus \mathbb{Z}$ of rank $2$, generated by two elements $x$ and $y$
subject to the single defining relation $[x,y] = 1$.
Then for every positive integer $m$, the subgroup $G_m=\lg x^m, y^m\rg$ is characteristic in $G$, with index $|G:G_m|=m^2$.
\end{prop}

Finally, we will use some Reidemeister-Schreier theory, which produces a defining presentation
for a subgroup $H$ of finite index in a finitely-presented group $G$.
An easily readable reference for this is~\cite[Chapter IV]{Johnson}, but in practice we use its implementation
as the {\tt Rewrite} command in the {\sc Magma} computation system~\cite{BCP97}.
We also found the groups that we use in the next section with the help of {\sc Magma} in constructing
and analysing some small examples.

\section{Main results}\label{Main results}

\begin{theorem}
\label{maintheorem}
For every  positive integer $m \geq 1$, there exist chiral $4$-polytopes $\P_m$ and $\Q_m$ of type $\{4,4,4\}$
with solvable automorphism groups of order $1024m^2$ and $2048m^2$, respectively.
\end{theorem}

\demo
We begin by defining ${\cal U}$ as the finitely-presented group \\[-10pt]
$$\lg\, a, b, c \ | \ a^4 = b^4 = c^4 = (ab)^2 = (bc)^2 = (abc)^2 = (a^2b^2)^4 = a^2c^2b^2(ac)^2 = [a,c^{-1}]b^2 = 1 \,\rg.$$
This group ${\cal U}$ has two normal subgroups of index $1024$ and $2048$, namely the subgroups
generated by $\{(ac^{-1})^4, (c^{-1}a)^4\}$ and $\{(bc^{-1})^4, (c^{-1}b)^4\}$, respectively.
The quotients of ${\cal U}$ by each of these give the initial members of our two infinite families.

\medskip
\noindent Case (1):  Take $N$ as the subgroup of ${\cal U}$ generated by $x = (ac^{-1})^4$ and $y = (c^{-1}a)^4$.

\smallskip
A short computation with {\sc Magma} shows that $N$ is normal in ${\cal U}$, with index $1024$. 
In fact, the defining relations for ${\cal U}$ can be used to show that  \\[-15pt]
\begin{center}
\begin{tabular}{lll}
$a^{-1}xa = y$, \qquad & $b^{-1}xb = y$ \quad &  and \quad $c^{-1}xc = y$,  \\[+2pt]
$a^{-1}ya = x$, \qquad & $b^{-1}yb = x^{-1}$ \quad & and \quad $c^{-1}yc = x^{-1}$.  \\[-3pt]
\end{tabular}
\end{center}

The first, third and fifth of these are easy to prove by hand, while the second and fourth can 
be verified in a number of ways, and the sixth follows from the other five. 
One way to prove the second and fourth is by hand, which we leave as a challenging exercise for the interested reader.
Another is by  a partial enumeration of cosets of the identity subgroup in ${\cal U}$. 
For example, if this is done using the {\tt ToddCoxeter} command in {\sc Magma}, allowing the definition of just 8000 cosets, 
then multiplication by each of the words $b^{-1}xby^{-1}$ and $a^{-1}yax^{-1}$ is found to fix the trivial coset, 
and therefore $b^{-1}xby^{-1} = 1 = a^{-1}yax^{-1}.$ 

It follows that conjugation by $a$, $b$ and $c$ induce the three permutations $\,(x,y)(x^{-1},y^{-1})$,
$\,(x,y,x^{-1},y^{-1})\,$ and $\,(x,y,x^{-1},y^{-1})$ on the set $\{x,y,x^{-1},y^{-1}\}$,
and then $(ac^{-1})^2$ and $(c^{-1}a)^2$ centralise both $x$ and $y$, so $x$ and $y$ centralise each other.

Also the {\tt Rewrite} command in {\sc Magma} gives a defining presentation for $N$, with $[x,y] = 1$ as a single 
defining relation.  Hence the normal subgroup $N$ is free abelian of rank 2. 


\smallskip
The quotient ${\cal U}/N$ is isomorphic to the automorphism group of the chiral 4-polytope
of type $\{4,4,4\}$ with 1024 automorphisms listed at~\cite{atles1}.

\smallskip
Now for any positive integer $m$, let $N_m$ be the subgroup generated by $x^m = (ac^{-1})^{4m}$
and $y^m = (c^{-1}a)^{4m}$. By Proposition~\ref{freeabelian}, we know that $N_m$ is characteristic in $N$
and hence normal in ${\cal U}$, with index $|\,{\cal U}:N_m| = |\,{\cal U}:N||N:N_m| = 1024m^2$.
Moreover, in the quotient $G_m = {\cal U}/N_m$, the subgroup  $N/N_m$ is abelian and normal,
with quotient $({\cal U}/N_m)/(N/N_m) \cong  {\cal U}/N$ being a $2$-group, and so $G_m$ is solvable,
by Proposition~\ref{solvable}.

Next, we use Proposition~\ref{quotient criterion} to prove that the triple $(\bar a, \bar b, \bar c)$ of
images of $a,b,c$ in $G_m$ satisfies the chiral form of the intersection condition.
To do this, we observe that the group with presentation $\lg \, u,v \ | \ u^4 = v^4 = (uv)^2 =(u^{2}v^{2})^4=1 \, \rg$
has order $2^7 = 128$, as does its image in the group $G_1 = {\cal U}/N_1$ under the epimorphism
taking $(u,v) \mapsto (aN_1,bN_1)$.  These claims are easily verifiable using {\sc Magma}.
Then since $G_1 = {\cal U}/N_1$ is a quotient of $G_m = {\cal U}/N_m$, the subgroup generated
by $\bar a$ and $\bar b$ in $G_m = {\cal U}/N_m$ must have order $128$ as well,
and hence the restriction to $\lg \bar a, \bar b \rg$ of the natural homomorphism from $G_m$ to $G_1$ is injective,
as required.

Accordingly, $G_m$ is the rotation group of a chiral or regular 4-polytope $\P_m$ of type $\{4,4,4\}$.
Assume for the moment that $\P_m$ is regular. Then there exists an automorphism $\rho$ of $G_m$
taking $(\bar a, \bar b, \bar c)$ to $(\bar a^{-1}, \bar a^{\,2}\bar b, \bar c)$, and  it  follows from the
relation $1 = \overline{a^2c^2b^2(ac)^2} = \bar a^2 \bar c^2 \bar b^2 (\bar a \bar c)^2$ that also
$1 = (\bar a^2 \bar c^2 \bar b^2 (\bar a \bar c)^2)^\rho =  \bar a^{-2} \bar c^{2} (\bar a^{2} \bar b)^2 (\bar a^{-1} \bar c)^2$ in $G_m$.
But the image of this element in $G_1$ has order 2, and as this is also the image of $a^{-2}c^{2}(a^{2}b)^{2}(a^{-1}c)^2)$ in $G_1$,
it follows that $(\bar a^2 \bar c^2 \bar b^2 (\bar a \bar c)^2)^\rho =  \bar a^{-2} \bar c^{2} (\bar a^{2} \bar b)^2 (\bar a^{-1} \bar c)^2$
is non-trivial in $G_m$, a contradiction. 

\smallskip
Thus $\P_m$ is chiral, with automorphism group $G_m$ of order $1024m^2$.

\medskip
\noindent Case (2):  Take $K$ as the subgroup of ${\cal U}$ generated by $z = (bc^{-1})^4$ and $w = (c^{-1}b)^4$.

\smallskip
Another computation with {\sc Magma} shows that $K$ is normal in ${\cal U}$, with index $2048$,
and moreover, the {\tt Rewrite} command tells us that $K$ is free abelian of rank $2$. 
In this case,  the defining relations for ${\cal U}$ give  \\[-18pt]
\begin{center}
\begin{tabular}{lll}
$a^{-1}za = z^{-1}$, \qquad & $b^{-1}zb = w$ \quad &  and \quad $c^{-1}zc = w$,  \\[+2pt]
$a^{-1}wa = w$, \qquad & $b^{-1}wb = z^{-1}$ \quad & and \quad $c^{-1}wc = z^{-1}$.  \\[-6pt]
\end{tabular}
\end{center}
The quotient ${\cal U}/K$ is the automorphism group of the chiral 4-polytope of type $\{4,4,4\}$
with 2048 automorphisms found by Zhang in~\cite{Zhang'sPhdThesis}.

\smallskip
Now for any positive integer $m$, let $K_m$ be the subgroup generated by $z^m = (bc^{-1})^{4m}$
and $w^m = (c^{-1}b)^{4m}$. Using Proposition~\ref{freeabelian}, we find that $K_m$ is characteristic in $K$
and hence normal in ${\cal U}$, with index $|\,{\cal U}:K_m| = |\,{\cal U}:K||K:K_m| = 2048m^2$.
Also the quotient $H_m = {\cal U}/K_m$ is solvable, again by Proposition~\ref{solvable}.

Next, the image of the subgroup generated by $a$ and $b$ in $H_1 = {\cal U}/K$ has order $128$,
so just as in Case (1) above, we can apply Proposition~\ref{quotient criterion} and find that the triple $(\bar a, \bar b, \bar c)$
of images of $a,b,c$ in $H_m = {\cal U}/K_m$ satisfies the chiral form of the intersection condition.

Thus $H_m$ is the rotation group of a chiral or regular 4-polytope $\Q_m$ of type $\{4,4,4\}$.

\smallskip
Moreover, the same argument as used in Case (1) shows that $\Q_m$ is chiral, because
the image in $H_1$ of the element $a^{-2}c^{2}(a^{2}b)^{2}(a^{-1}c)^2)$ has order $2$,
and hence the image of that element in $H_m$ is non-trivial.

\smallskip
Thus $\Q_m$ is chiral, with automorphism group $H_m$ of order $2048m^2$.
 \hfill\qed

As a special case we have the following Corollary, which is an immediate consequence of Theorem~\ref{maintheorem}
when $m$ is taken as a power of $2$.

\begin{cor}
There exists a chiral $4$-polytope of type $\{4,4,4\}$ with automorphism group of order $2^n,$
for every integer $n \ge 10$.
\end{cor}

On the other hand, an inspection of the lists at \cite{atles1} shows that there exists no 
such chiral polytope with automorphism group of order $2^n$, where $n \le 9$.  

\newpage\f
%
{\bf\Large Acknowledgements}

\medskip\smallskip
\noindent
The first author acknowledges the hospitality of Beijing Jiaotong University,
and partial support from the  N.Z. Marsden Fund (project UOA1626). 
The second and the third authors acknowledge the partial support from 
the National Natural Science Foundation of China (11731002) and the 111 Project of China (B16002). The three authors also acknowledge the extensive use of the {\sc Magma} system~\cite{BCP97} in helping conduct experiments that helped show the way towards the main theorem.
\\[-12pt]



\begin{thebibliography}{99}

\bibitem{BCP97}
W. Bosma, J. Cannon and C. Playoust:
The {M}agma {A}lgebra {S}ystem. {I}:  the user language.
{\em J. Symbolic Comput.} 24 (1997) 235--265.
\\[-20pt]

\bibitem{ChiralPolytopesofPSLextension}
A. Breda D'Azevedo, G.A. Jones, E. Schulte:
Constructions of chiral polytopes of small rank,
{\em Canad. J. Math.} 63 (2011) 1254--1283.
\\[-20pt]

\bibitem{ConstructChiralPolytopes}
M.D.E. Conder, I. Hubard, T. Pisanski:
Constructions for chiral polytopes,
{\em J. London Math. Soc.} 77 (2008) 115--129.
\\[-20pt]

\bibitem{atles1}
M.D.E. Conder:
Chiral polytopes with up to 4000 flags,
\url{https://www.math.auckland.ac.nz/~conder/ChiralPolytopesWithUpTo4000Flags-ByOrder.txt}.
\\[-20pt]

\bibitem{ChiralPolytopesOfANSN}
M.D.E. Conder, I. Hubard, E. O'Reilly Regueiro, D. Pellicer:
Construction of chiral 4-polytopes with an alternating or symmetric group as automorphism group,
{\em J. Algebraic Combin.} 42 (2015) 225--244.
\\[-20pt]

\bibitem{AbelianCoverOfChiralPolytopes}
M.D.E. Conder, W.-J. Zhang:
Abelian covers of chiral polytopes,
{\em J. Algebra\/} 478 (2017) 437--457.
\\[-20pt]

\bibitem{BooksCoxeter}
H.S.M. Coxeter, W.O.J. Moser:
{\em Generators and Relations for Discrete Groups\/}, 4th edition,
Springer, Berlin, 1980.
\\[-20pt]

\bibitem{TightChiralPolyhedra}
G. Cunningham:
Tight chiral polyhedra,
{\em Combinatorica\/} 38 (2018) 115--142.
\\[-20pt]

\bibitem{moreonConstructChiralPolytopesforanyrank}
G. Cunningham, D. Pellicer:
Chiral extensions of chiral polytopes,
{\em Discrete Math.} 330 (2014) 51--60.
\\[-20pt]

\bibitem{atles2}
M. Hartley, I. Hubard, D. Leemans:
An Atlas of Chiral Polytopes for Small Almost Simple Groups,
\url{http://homepages.ulb.ac.be/~dleemans/CHIRAL/index.html}.
\\[-20pt]

\bibitem{HFL}
D.-D. Hou, Y.-Q. Feng, D. Leemans:
Existence of regular $3$-polytopes of order $2^n$,
{\em J. Group Theory\/} 22 (2019) 579--616.
\\[-20pt]

\bibitem{HFL1}
D.-D. Hou, Y.-Q. Feng, D. Leemans:
On regular polytopes of $2$-powers,
{\em Discrete Comput. Geom.\/}, to appear; arXiv:1901.07200.
\\[-20pt]

\bibitem{Johnson}
D.L. Johnson,
{\em Topics in the Theory of Group Presentations},
Cambridge Univ.~Press, Cambridge (1980).
\\[-20pt]

\bibitem{ARP}
P. McMullen, E. Schulte:
{\em Abstract Regular Polytopes\/}, Encyclopedia Math. Appl., vol. 92,
Cambridge University Press, Cambridge, 2002.
\\[-20pt]

\bibitem{BE2009}
B. Monson, E. Schulte:
Modular reduction in abstract polytopes,
{\em Canad. Math. Bull.} 52 (2009) 435--450.
\\[-20pt]

\bibitem{ConstructChiralPolytopesforanyrank}
D. Pellicer:
A construction of higher rank chiral polytopes,
{\em Discrete Math.} 310 (2010) 1222--1237.
\\[-20pt]

\bibitem{ChiralPolytopes}
E. Schulte, A.I. Weiss:
Chiral polytopes, in:
{\em Applied Geometry and Discrete Mathematics\/},
DI-MACS Ser. Discrete Math. Theoret. Comput. Sci., vol.4,
Amer. Math. Soc., Providence, RI, (1991) 493--516.
\\[-20pt]

\bibitem{ChiralPolytopesWithPSL}
E. Schulte, A.I. Weiss:
Chirality and projective linear groups,
{\em Discrete Math.} 131 (1994) 221--261.
\\[-20pt]

\bibitem{FEChiralPolytopes}
E. Schulte, A.I. Weiss:
Free extensions of chiral polytopes,
{\em Canad. J. Math.} 47 (1995) 641--654.
\\[-20pt]

\bibitem{Problem}
E. Schulte,  A.I. Weiss:
Problems on polytopes, their groups, and realizations,
{\em Periodica Math. Hungarica\/} 53 (2006) 231--255.
\\[-20pt]

\bibitem{GroupBook}
M.Y. Xu:
{\em Introduction to Group Theory I\/},
Science Publishing House, Beijing, 1999.
\\[-20pt]

\bibitem{Zhang'sPhdThesis}
W.-J. Zhang:
{\em Constructions For Chiral Polytopes\/},
PhD Thesis, University of Auckland, 2016.

\end{thebibliography}
\end{document}